\documentclass{article}%
\usepackage{amsfonts}
\usepackage{amsmath}
\usepackage{amssymb}
\usepackage{graphicx}%
\providecommand{\U}[1]{\protect\rule{.1in}{.1in}}
\newtheorem{theorem}{Theorem}

\newtheorem{lemma}[theorem]{Lemma}

\newenvironment{proof}[1][Proof]{\noindent\textbf{#1.} }{\ \rule{0.5em}{0.5em}}
\begin{document}

\title{A Spectral Method for Elliptic Equations: \\The Dirichlet Problem}
\author{Kendall Atkinson\\Departments of Mathematics \& Computer Science \\The University of Iowa
\and David Chien, Olaf Hansen\\Department of Mathematics \\California State University San Marcos}
\maketitle

\begin{abstract}
Let $\Omega$ be an open, simply connected, and bounded region in
$\mathbb{R}^{d}$, $d\geq2$, and assume its boundary $\partial\Omega$ is
smooth. Consider solving an elliptic partial differential equation $Lu=f$ over
$\Omega$ with zero Dirichlet boundary values. The problem is converted to an
equivalent\ elliptic problem over the unit ball $B$, and then a spectral
method is given that uses a special polynomial basis. With sufficiently smooth
problem parameters, the method is shown to have very fast convergence.
Numerical examples illustrate exponential convergence.

\end{abstract}

\section{INTRODUCTION}

Consider\ solving the elliptic partial differential equation%
\begin{equation}
Lu(\mathbf{s})\equiv-\sum_{i,j=1}^{d}\frac{\partial}{\partial s_{i}}\left(
a_{i,j}(\mathbf{s})\frac{\partial u(\mathbf{s})}{\partial s_{j}}\right)
+\gamma(\mathbf{s})u(\mathbf{s})=f(\mathbf{s}),\quad\quad\mathbf{s}%
\in\text{$\Omega$}\subseteq\mathbb{R}^{d} \label{e1}%
\end{equation}
with the Dirichlet boundary condition
\begin{equation}
u(\mathbf{s})\equiv0,\quad\quad\mathbf{s}\in\partial\Omega\label{e2}%
\end{equation}
Assume $d\geq2$. Assume $\Omega$ is an open, simply-connected, and bounded
region in $\mathbb{R}^{d}$, and assume that its boundary $\partial\Omega$ is
several times continuously differentiable. Similarly, assume the functions
$\gamma(\mathbf{s}),$ $f(\mathbf{s}),$ $a_{i,j}(\mathbf{s})$ are several times
continuously differentiable over $\overline{\Omega}$. As usual, assume the
matrix $A(\mathbf{s})=\left[  a_{i,j}(\mathbf{s})\right]  $ is symmetric, and
also assume it satisfies the strong ellipticity condition,%
\begin{equation}
\xi^{\text{T}}A(\mathbf{s})\xi\geq c_{0}\xi^{\text{T}}\xi,\quad\quad
\mathbf{s}\in\overline{\Omega},\quad\xi\in\mathbb{R}^{d} \label{e3}%
\end{equation}
with $c_{0}>0$. Also assume $\gamma(\mathbf{s})\geq0$, $\mathbf{s}\in\Omega$.

In \S \ref{spec_method} we consider the special region $\Omega=B$, the open
unit ball in $\mathbb{R}^{d}$. We define a Galerkin method for (\ref{e1}%
)-(\ref{e2}) with a special finite-dimensional subspace of polynomials, and we
give an error analysis that shows rapid convergence of the method. In
\S \ref{elliptic_transf} we discuss the use of a transformation from a general
region $\Omega$ to the unit ball $B$, showing that the transformed equation is
again elliptic over $B$. \ Implementation issues are discussed in
\S \ref{implementation} for problems in $\mathbb{R}^{2}$ and $\mathbb{R}^{3}$.
\ We conclude in \S \ref{num_example} with numerical examples in
$\mathbb{R}^{2}$ and $\mathbb{R}^{3}$.

The methods of this paper generalize to the equation%
\begin{align*}
Lu(\mathbf{s})\equiv &  -\sum_{i,j=1}^{d}\frac{\partial}{\partial s_{i}%
}\left(  a_{i,j}(\mathbf{s})\frac{\partial u(\mathbf{s})}{\partial s_{j}%
}\right) \\
&  +\sum_{j=1}^{d}b_{j}\left(  \mathbf{s}\right)  \frac{\partial
u(\mathbf{s})}{\partial s_{j}}+\gamma(\mathbf{s})u(\mathbf{s})=f(\mathbf{s}%
),\quad\quad\mathbf{s}\in\text{$\Omega$}\subseteq\mathbb{R}^{d}%
\end{align*}
which contains first order derivative terms, provided the operator $L$ is
strongly elliptic. To do so, use the results of Brenner and Scott
\cite[\S \S 2.6-2.8]{brenner-scott}, combined with the methods of the present
paper. We have chosen to restrict our work to the more standard symmetric
problem (\ref{e1}).

There is a rich literature on spectral methods for solving partial
differential equations. From the more recent literature, we cite \cite{doha},
\cite{shen99}, and \cite{shen07}. Their bibliographies contain \ references to
earlier papers on spectral methods. Our approach is somewhat different than
the standard approaches, as we are converting the partial differential
equation to an equivalent problem on the unit disk or unit ball, and in the
process we are required to work with a more complicated equation. Our approach
is reminiscent of the use of conformal mappings for planar problems. Conformal
mappings can be used with our approach when working on planar problems,
although having a conformal mapping is not necessary.

\section{A spectral method on the unit ball\label{spec_method}}

The Dirichlet problem (\ref{e1})-(\ref{e2}) has the following variational
reformulation: Find $u\in H_{0}^{1}\left(  \Omega\right)  $ such that%
\begin{equation}%
\begin{array}
[c]{r}%
{\displaystyle\int\nolimits_{\Omega}}
\left[
{\displaystyle\sum\limits_{i,j=1}^{d}}
a_{i,j}(\mathbf{s})\dfrac{\partial u(\mathbf{s})}{\partial s_{j}}%
\dfrac{\partial v(\mathbf{s})}{\partial s_{i}}+\gamma(\mathbf{s}%
)u(\mathbf{s})v(\mathbf{s})\right]  d\mathbf{s}\quad\quad\quad\quad\medskip\\
=%
{\displaystyle\int_{\Omega}}
f(\mathbf{s})v(\mathbf{s})\,d\mathbf{s},\quad\quad\forall v\in H_{0}%
^{1}\left(  \Omega\right)
\end{array}
\label{e20}%
\end{equation}
We define a spectral Galerkin method in this section for the special region
$\Omega=B$. In \S \ref{elliptic_transf} we discuss the transformation of
(\ref{e1}) from a general $\overline{\Omega}$ to an equivalent equation over
the unit ball $\overline{B}$, a transformation that retains the ellipticity of
the problem. In the remainder of this section, we replace $\Omega$ with $B$.

Introduce the bilinear form%
\begin{equation}
\mathcal{A}\left(  v,w\right)  =\int_{B}\left[  \sum_{i,j=1}^{d}%
a_{i,j}(\mathbf{x})\frac{\partial v(\mathbf{x})}{\partial x_{j}}\frac{\partial
w(\mathbf{x})}{\partial x_{i}}+\gamma(\mathbf{x})v(\mathbf{x})w(\mathbf{x}%
)\right]  d\mathbf{x},\quad v,w\in H_{0}^{1}\left(  B\right)  \label{e21}%
\end{equation}
and the bounded linear functional%
\[
\ell(v)=\int_{B}f\left(  \mathbf{x}\right)  v\left(  \mathbf{x}\right)
\,d\mathbf{x},\quad\quad v\in H_{0}^{1}\left(  B\right)
\]
The variational problem (\ref{e20}) can now be written as follows: find $u\in
H_{0}^{1}\left(  B\right)  $ for which%
\begin{equation}
\mathcal{A}\left(  u,v\right)  =\ell(v),\quad\quad\forall v\in H_{0}%
^{1}\left(  B\right)  \label{e22}%
\end{equation}

It is straightforward to show $\mathcal{A}$ is bounded,%
\[%
\begin{array}
[c]{c}%
\left\vert \mathcal{A}\left(  v,w\right)  \right\vert \leq c_{\mathcal{A}%
}\left\Vert v\right\Vert _{1}\left\Vert w\right\Vert _{1}\medskip\\
c_{\mathcal{A}}=\max\limits_{\mathbf{x}\in\overline{B}}\left\Vert
A(\mathbf{x})\right\Vert _{2}+\left\Vert \gamma\right\Vert _{\infty}%
\end{array}
\]
with $\left\Vert \cdot\right\Vert _{1}$ the norm of $H_{0}^{1}\left(
\Omega\right)  $ and $\left\Vert A(\mathbf{x})\right\Vert _{2}$ the matrix
2-norm of the matrix $A\left(  \mathbf{x}\right)  $.\ In addition, we assume
\begin{equation}
\mathcal{A}\left(  v,v\right)  \geq c_{e}\Vert v\Vert_{1}^{2},\quad\quad v\in
H_{0}^{1}\left(  B\right)  \label{e23}%
\end{equation}
This follows generally from (\ref{e3}) and the size of the function
$\gamma(\mathbf{x})$ over $\overline{B}$; when $\gamma\equiv0$, $c_{e}=c_{0}$.
Under standard assumptions on $\mathcal{A}$, including the strong ellipticity
in (\ref{e23}), the Lax-Milgram Theorem implies the existence of a unique
solution $u$ to (\ref{e22}) with
\[
\Vert u\Vert_{1}\leq\frac{1}{c_{e}}\Vert\ell\Vert
\]

Denote by $\Pi_{n}$ the space of polynomials in $d$ variables that are of
degree $\leq n$: $p\in\Pi_{n}$ if it has the form%
\[
p(\mathbf{x})=\sum_{\left\vert i\right\vert \leq n}a_{i}x_{1}^{i_{1}}%
x_{2}^{i_{2}}\dots x_{d}^{i_{d}}
\]
with $i$ a multi-integer, $i=\left(  i_{1},\dots,i_{d}\right)  $, and
$\left\vert i\right\vert =i_{1}+\cdots+i_{d}$. Let $\mathcal{X}_{n}$ denote
our approximation subspace,%
\begin{equation}
\mathcal{X}_{n}=\left\{  \left(  1-\left\Vert \mathbf{x}\right\Vert _{2}%
^{2}\right)  p(\mathbf{x})\mid p\in\Pi_{n}\right\}  \label{e26}%
\end{equation}
with $\left\Vert \mathbf{x}\right\Vert _{2}^{2}=x_{1}^{2}+\cdots+x_{d}^{2}$.
The subspaces $\Pi_{n}$ and $\mathcal{X}_{n}$ have dimension%
\[
N_{n}=\binom{n+d}{d}
\]

\begin{lemma}
\label{lemma1}Let $\Delta$ denote the Laplacian operator in $\mathbb{R}^{d}$.
Then%
\begin{equation}
\Delta:\mathcal{X}_{n}\overset{1-1}{\underset{onto}{\longrightarrow}}\Pi_{n}
\label{e28}%
\end{equation}

\end{lemma}

\noindent For a short proof, see \cite{atkinson-hansen}.

The Galerkin method for obtaining an approximate solution to (\ref{e22}) is as
follows: find $u_{n}\in\mathcal{X}_{n}$ for which%
\begin{equation}
\mathcal{A}\left(  u_{n},v\right)  =\ell(v),\quad\quad\forall v\in
\mathcal{X}_{n} \label{e30}%
\end{equation}
The Lax-Milgram Theorem (cf. \cite[\S 8.3]{atkinson-han}, \cite[\S 2.7]%
{brenner-scott}) implies the existence of $u_{n}$ for all $n$. For the error
in this Galerkin method, Cea's Lemma (cf. \cite[p. 365]{atkinson-han},
\cite[p. 62]{brenner-scott}) implies the convergence of $u_{n}$ to $u$, and
moreover,%
\begin{equation}
\Vert u-u_{n}\Vert_{1}\leq\frac{c_{\mathcal{A}}}{c_{e}}\inf_{v\in
\mathcal{X}_{n}}\Vert u-v\Vert_{1} \label{e31}%
\end{equation}
It remains to bound the best approximation error on the right side of this inequality.

Given an arbitrary $u\in H_{0}^{2}\left(  B\right)  $, define $w=-\Delta u$.
Then $w\in L^{2}\left(  B\right)  $ and $u$ satisfies the boundary value
problem
\begin{align*}
-\Delta u(P)  &  =w(P),\quad\quad P\in\text{$B$}\\
u(P)  &  =0,\quad\quad P\in\partial B
\end{align*}
It follows that
\begin{equation}
u(P)=\int_{B}G(P,Q)w(Q)\,dQ,\quad\quad P\in\overline{B} \label{e41}%
\end{equation}
For $\mathbb{R}^{2}$ and $\mathbb{R}^{3}$, the Green's function is defined as
follows.%
\begin{equation}%
\begin{array}
[c]{cl}%
d=2: & G(P,Q)=\dfrac{1}{2\pi}\log\dfrac{\left\vert P-Q\right\vert }{\left\vert
\mathcal{T}(P)-Q\right\vert },\quad\quad\medskip\\
d=3: & G(P,Q)=-\dfrac{1}{4\pi}\left\{  \dfrac{1}{\left\vert P-Q\right\vert
}-\dfrac{1}{\left\vert P\right\vert }\dfrac{1}{\left\vert \mathcal{T}%
(P)-Q\right\vert }\right\}
\end{array}
\label{e42}%
\end{equation}
for $P\neq Q,$ $Q\in B,$ $P\in\overline{B}$. $\mathcal{T}(P)$ denotes the
inverse point for $P$ with respect to the unit sphere $S^{d-1}\subseteq
\mathbb{R}^{d}$,%
\[
\mathcal{T}(r\mathbf{x})=\frac{1}{r}\mathbf{x},\quad\quad0<r\leq
1,\quad\mathbf{x}\in S^{d-1}
\]
Differentiate (\ref{e41}) to obtain%
\begin{equation}
\nabla u\left(  P\right)  =\int_{B}\left[  \nabla_{P}G(P,Q)\right]
w(Q)\,dQ,\quad\quad P\in\overline{B} \label{e43}%
\end{equation}
Note that $\nabla_{P}G(P,\cdot)$ is absolutely integrable over $\overline{B}$,
for all $P\in\overline{B}$.

Let $w_{n}\in\Pi_{n}\,$\ be an approximation of $w$, say in the norm of either
$C\left(  \overline{B}\right)  $ or $L^{2}\left(  B\right)  $, and let
\[
q_{n}(P)=\int_{B}G(P,Q)w_{n}(Q)\,dQ,\quad\quad P\in\overline{B}
\]
We can show $q_{n}\in\mathcal{X}_{n}$. This follows from Lemma \ref{lemma1}
and noting that the mapping in (\ref{e41}) is the inverse of (\ref{e28}).

Then we have
\[
u(P)-q_{n}(P)=\int_{B}G(P,Q)\left[  w(P)-w_{n}(Q)\right]  \,dQ,\quad\quad
P\in\overline{B}
\]%
\[
\nabla\left[  u\left(  P\right)  -q_{n}(P)\right]  =\int_{B}\left[  \nabla
_{P}G(P,Q)\right]  \left[  w(Q)-w_{n}(Q)\right]  \,dQ,\quad\quad P\in
\overline{B}
\]
The integral operators on the right side are weakly singular compact integral
operators on $L^{2}\left(  B\right)  $ to $L^{2}\left(  B\right)  $
\cite[Chap. 7, \S 3]{mikhlin}. This implies%
\begin{equation}
\Vert u-q_{n}\Vert_{1}\leq c\Vert w-w_{n}\Vert_{0} \label{e47}%
\end{equation}
By letting $w_{n}$ be the orthogonal projection of $w$ into $\Pi_{n}$, the
right side will go to zero since the polynomials are dense in $L^{2}\left(
B\right)  $. In turn, this implies convergence in the $H_{0}^{1}\left(
B\right)  $ norm for the right side in (\ref{e31}) provided $u\in H_{0}%
^{2}\left(  B\right)  $.

The result
\[
\inf_{v\in\mathcal{X}_{n}}\Vert u-v\Vert_{1}\rightarrow0\quad\,\text{as}\quad
n\rightarrow\infty,\quad\quad u\in H_{0}^{2}\left(  B\right)
\]
can be extended to any $u\in H_{0}^{1}\left(  B\right)  $. It basically
follows from the denseness of $H_{0}^{2}\left(  B\right)  $ in $H_{0}%
^{1}\left(  B\right)  $. Let $u\in H_{0}^{1}\left(  B\right)  $. We need to
find a sequence of polynomials $\left\{  q_{n}\right\}  $ for which $\Vert
u-q_{n}\Vert_{1}\rightarrow0$. We know $H_{0}^{2}\left(  B\right)  $ is dense
in $H_{0}^{1}\left(  B\right)  $. Given any $k>0$, choose $u_{k}\in H_{0}%
^{2}\left(  B\right)  $ with $\Vert u-u_{k}\Vert_{1}\leq1/k$. Then choose a
polynomial $w_{k}$ for which we have the corresponding polynomial $q_{k}$
satisfying $\Vert u_{k}-q_{k}\Vert_{1}\leq1/k$, based on (\ref{e47}).
[Regarding the earlier notation, $q_{k}$ need not be of degree $\leq k$.] Then
$\Vert u-q_{k}\Vert_{1}\leq2/k$.

To obtain orders of convergence, use (\ref{e47}) and results on best
multivariate polynomial approximation over the unit disk. For example, use
results of Ragozin \cite[Thm 3.4]{ragozin} or Yuan Xu \cite{xu2006}. From
\cite{ragozin} we have the following.

\begin{theorem}
\label{theorem1}Assume $u\in C^{k+2}\left(  \overline{B}\right)  $ for some
$k>0$. Then there is a polynomial $q_{n}\in\mathcal{X}_{n}$ for which%
\begin{equation}
\left\Vert u-q_{n}\right\Vert _{\infty}\leq D\left(  k,d\right)  n^{-k}\left(
n^{-1}\left\Vert u\right\Vert _{\infty,k+2}+\omega\left(  u^{(k+2)}%
,1/n\right)  \right)  \label{e48}%
\end{equation}

\end{theorem}

\noindent In this,
\[
\left\Vert u\right\Vert _{\infty,k+2}=\sum_{\left\vert i\right\vert \leq
k+2}\left\Vert \partial^{i}u\right\Vert _{\infty}
\]%
\[
\omega\left(  g,\delta\right)  =\sup_{\left\vert \mathbf{x}-\mathbf{y}%
\right\vert \leq\delta}\left\vert g\left(  \mathbf{x}\right)  -g\left(
\mathbf{y}\right)  \right\vert
\]%
\[
\omega\left(  u^{(k+2)},\delta\right)  =\sum_{\left\vert i\right\vert
=k+2}\omega\left(  \partial^{i}u,\delta\right)
\]

\section{Transformation of the elliptic equation \label{elliptic_transf}}

Consider the differential operator%
\begin{equation}
\mathcal{M}v(\mathbf{s})=-\sum_{i,j=1}^{d}\frac{\partial}{\partial s_{i}%
}\left(  a_{i,j}(\mathbf{s})\frac{\partial v(\mathbf{s})}{\partial s_{j}%
}\right)  ,\quad\quad\mathbf{s}\in\text{$\Omega$}\subseteq\mathbb{R}^{d},\quad
v\in C^{2}\left(  \overline{\Omega}\right)  \label{e60}%
\end{equation}
which satisfies the ellipticity condition (\ref{e3}) with $c_{0}>0$. \ The
operator $\mathcal{M}$ is said to be elliptic on $H^{2}\left(  \Omega\right)
$.\ We want to transform the operator $\mathcal{M}$ to one acting on functions
$\widetilde{u}\in C^{2}(\overline{B})$ with $B$ the unit ball in
$\mathbb{R}^{d}$.

Assume the existence of a twice-differentiable mapping%
\begin{equation}
\Phi:\overline{B}\underset{onto}{\overset{1-1}{\longrightarrow}}%
\overline{\Omega} \label{e61}%
\end{equation}
and let $\Psi=\Phi^{-1}:\overline{\Omega}\underset{onto}{\overset
{1-1}{\longrightarrow}}\overline{B}$. Let%
\[
J(\mathbf{x})\equiv\left(  D\Phi\right)  (\mathbf{x})=\left[  \frac
{\partial\varphi_{i}(\mathbf{x})}{\partial x_{j}}\right]  _{i,j=1}^{d}%
,\quad\quad\mathbf{x}\in\overline{B}\subseteq\mathbb{R}^{d}%
\]
denote the Jacobian of the transformation. As usual we assume $J(\mathbf{x})$
is nonsingular on $\overline{B}$, and furthermore%
\begin{equation}
\min_{\mathbf{x}\in\overline{B}}\,\left\vert \det J(\mathbf{x})\right\vert >0
\label{e64}%
\end{equation}
Similarly, let $K(\mathbf{s})\equiv\left(  D\Psi\right)  (\mathbf{s})$ denote
the Jacobian of $\Psi$ over $\overline{\Omega}$. \ By differentiating the
components of the equation%
\[
\Psi\left(  \Phi(\mathbf{x})\right)  =\mathbf{x}%
\]
we obtain%
\[
K\left(  \Phi(\mathbf{x})\right)  =J^{-1}\left(  \mathbf{x}\right)
,\quad\quad\mathbf{x}\in\overline{B}%
\]
This general approach is reminiscent of the coordinate transformations in
\cite[Chap. 2]{liseikin} in which the mapping function is used in generating a
mesh on a region $\Omega$.

For $v\in C^{2}(\overline{\Omega})$, let%
\[
\widetilde{v}(\mathbf{x})=v\left(  \Phi\left(  \mathbf{x}\right)  \right)
,\quad\quad\mathbf{x}\in\overline{B}\subseteq\mathbb{R}^{d}
\]
and conversely,%
\[
v(\mathbf{s})=\widetilde{v}\left(  \Psi\left(  \mathbf{s}\right)  \right)
,\quad\quad\mathbf{s}\in\overline{\Omega}\subseteq\mathbb{R}^{d}
\]
Then
\begin{align*}
\frac{\partial\widetilde{v}}{\partial x_{i}}  &  =\frac{\partial v}{\partial
s_{1}}\frac{\partial\varphi_{1}\left(  \mathbf{x}\right)  }{\partial x_{i}%
}+\cdots+\frac{\partial v}{\partial s_{d}}\frac{\partial\varphi_{d}\left(
\mathbf{x}\right)  }{\partial x_{i}}\medskip\\
&  =\left[  \frac{\partial\varphi_{1}\left(  \mathbf{x}\right)  }{\partial
x_{i}},\cdots,\frac{\partial\varphi_{d}\left(  \mathbf{x}\right)  }{\partial
x_{i}}\right]  \nabla_{\mathbf{s}}v
\end{align*}
with the gradient $\nabla_{\mathbf{s}}v$ a column vector evaluated at
$\mathbf{s}=\Phi\left(  \mathbf{x}\right)  $. More concisely,%
\begin{equation}
\nabla_{\mathbf{x}}\widetilde{v}\left(  \mathbf{x}\right)  =J\left(
\mathbf{x}\right)  ^{\text{T}}\nabla_{\mathbf{s}}v\left(  \mathbf{s}\right)
,\quad\quad\mathbf{s}=\Phi\left(  \mathbf{x}\right)  \label{e66}%
\end{equation}
Similarly,%
\begin{equation}
\nabla_{\mathbf{s}}v(\mathbf{s})=K(\mathbf{s})^{\text{T}}\nabla_{\mathbf{x}%
}\widetilde{v}(\mathbf{x}),\quad\quad\mathbf{x}=\Psi(\mathbf{s}) \label{e67}%
\end{equation}

\begin{theorem}
\label{theorem2}Assume the transformation $\Phi$ satisfies (\ref{e61}) or
(\ref{e64}). Then for $\mathbf{s}=\Phi(\mathbf{x})$,%
\begin{equation}
\left(  \mathcal{M}v\right)  (\mathbf{s})=-\frac{1}{\det\left(  J(\mathbf{x}%
)\right)  }\sum_{i,j=1}^{d}\frac{\partial}{\partial x_{i}}\left(  \det\left(
J(\mathbf{x})\right)  \widetilde{a}_{i,j}(\mathbf{x})\frac{\partial
\widetilde{v}(\mathbf{x})}{\partial x_{j}}\right)  \label{e68}%
\end{equation}%
\begin{align}
\widetilde{A}\left(  \mathbf{x}\right)   &  =K\left(  \Phi\left(
\mathbf{x}\right)  \right)  A(\Phi\left(  \mathbf{x}\right)  )K\left(
\Phi\left(  \mathbf{x}\right)  \right)  ^{\text{T}}\label{e69}\\
&  \equiv\left[  \widetilde{a}_{i,j}(\mathbf{x})\right]  _{i,j=1}^{d}\nonumber
\end{align}

\begin{proof}
Let $w\in C_{0}^{\infty}\left(  \overline{\Omega}\right)  $. Then
\begin{equation}
\int_{\Omega}\left(  \mathcal{M}v\right)  (\mathbf{s})w(\mathbf{s}%
)\,d\mathbf{s}=\int_{B}\left(  \mathcal{M}v\right)  (\Phi\left(
\mathbf{x}\right)  )w(\Phi\left(  \mathbf{x}\right)  )\det\left(
J(\mathbf{x})\right)  \,d\mathbf{x} \label{e70}%
\end{equation}
On the other hand, using integration by parts we have%
\begin{align}
\int_{\Omega}\left(  \mathcal{M}v\right)  (\mathbf{s})w(\mathbf{s}%
)\,d\mathbf{s}  &  =\int_{\Omega}\sum_{i,j=1}^{d}a_{i,j}(\mathbf{s}%
)\frac{\partial v(\mathbf{s})}{\partial s_{j}}\frac{\partial w(\mathbf{s}%
)}{\partial s_{i}}\,d\mathbf{s}\nonumber\\
&  =\int_{B}\sum_{i,j=1}^{d}a_{i,j}(\Phi\left(  \mathbf{x}\right)
)\frac{\partial v(\Phi\left(  \mathbf{x}\right)  )}{\partial s_{j}}%
\frac{\partial w(\Phi\left(  \mathbf{x}\right)  )}{\partial s_{i}}\det\left(
J(\mathbf{x})\right)  \,d\mathbf{x} \label{e71}%
\end{align}
Using (\ref{e67}),%
\[%
\begin{array}
[c]{l}%
{\displaystyle\sum\limits_{i,j=1}^{d}}
a_{i,j}(\Phi\left(  \mathbf{x}\right)  )\dfrac{\partial v(\Phi\left(
\mathbf{x}\right)  )}{\partial s_{j}}\dfrac{\partial w(\Phi\left(
\mathbf{x}\right)  )}{\partial s_{i}}=\left[  \nabla_{\mathbf{s}}w(\Phi\left(
\mathbf{x}\right)  )\right]  ^{\text{T}}A\left(  \Phi\left(  \mathbf{x}%
\right)  \right)  \left[  \nabla_{\mathbf{s}}v(\Phi\left(  \mathbf{x}\right)
)\right]  \medskip\\
\quad\quad\quad\quad\quad\quad=\left[  \nabla_{\mathbf{x}}\widetilde
{w}(\mathbf{x})\right]  ^{\text{T}}K(\Phi\left(  \mathbf{x}\right)  )A\left(
\Phi\left(  \mathbf{x}\right)  \right)  K(\Phi\left(  \mathbf{x}\right)
)^{\text{T}}\left[  \nabla_{\mathbf{x}}\widetilde{v}(\mathbf{x})\right]
\medskip\\
\quad\quad\quad\quad\quad\quad=\left[  \nabla_{\mathbf{x}}\widetilde
{w}(\mathbf{x})\right]  ^{\text{T}}\widetilde{A}\left(  \mathbf{x}\right)
\left[  \nabla_{\mathbf{x}}\widetilde{v}(\mathbf{x})\right]
\end{array}
\]
Using this to continue (\ref{e71}),%
\begin{align}
\int_{\Omega}\left(  \mathcal{M}v\right)  (\mathbf{s})w(\mathbf{s}%
)\,d\mathbf{s}  &  =\int_{B}\left[  \nabla_{\mathbf{x}}\widetilde
{w}(\mathbf{x})\right]  ^{\text{T}}\widetilde{A}\left(  \mathbf{x}\right)
\left[  \nabla_{\mathbf{x}}\widetilde{v}(\mathbf{x})\right]  \det\left(
J(\mathbf{x})\right)  \,d\mathbf{x}\nonumber\\
&  =\int_{\Omega}\sum_{i,j=1}^{d}\widetilde{a}_{i,j}(\mathbf{s})\frac
{\partial\widetilde{v}(\mathbf{x})}{\partial x_{j}}\frac{\partial\widetilde
{w}(\mathbf{x})}{\partial x_{i}}\det\left(  J(\mathbf{x})\right)
\,d\mathbf{x}\nonumber\\
&  =-\int_{\Omega}\sum_{i,j=1}^{d}\frac{\partial}{\partial x_{i}}\left(
\det\left(  J(\mathbf{x})\right)  \widetilde{a}_{i,j}(\mathbf{x}%
)\frac{\partial\widetilde{v}(\mathbf{x})}{\partial x_{j}}\right)
\widetilde{w}(\mathbf{x})\,d\mathbf{x} \label{e73}%
\end{align}
Comparing (\ref{e70}) and (\ref{e73}), and noting that $w\in C_{0}^{\infty
}\left(  \overline{\Omega}\right)  $ is arbitrary, we have%
\[
\left(  \mathcal{M}v\right)  (\Phi\left(  \mathbf{x}\right)  )\det\left(
J(\mathbf{x})\right)  =-\sum_{i,j=1}^{d}\frac{\partial}{\partial x_{i}}\left(
\det\left(  J(\mathbf{x})\right)  \widetilde{a}_{i,j}(\mathbf{x}%
)\frac{\partial\widetilde{v}(\mathbf{x})}{\partial x_{j}}\right)
\]
which proves (\ref{e68}). \hfill
\end{proof}
\end{theorem}

With this transformation, we can solve the Dirichlet problem over a general
region $\Omega$ by transforming it to an equivalent problem over the unit ball
$B$. We can apply the Galerkin method to (\ref{e1}) by means of the
transformation (\ref{e68}). We convert (\ref{e1}) to the equation%
\begin{equation}%
\begin{array}
[c]{r}%
-%
{\displaystyle\sum\limits_{i,j=1}^{d}}
\dfrac{\partial}{\partial x_{i}}\left(  \det\left(  J(\mathbf{x})\right)
\widetilde{a}_{i,j}(\mathbf{x})\dfrac{\partial\widetilde{v}(\mathbf{x}%
)}{\partial x_{j}}\right)  +\det\left(  J(\mathbf{x})\right)  \gamma
(\Phi\left(  \mathbf{x}\right)  )\widetilde{v}(\mathbf{x})\quad\quad\quad\\
=\det\left(  J(\mathbf{x})\right)  f\left(  \Phi\left(  \mathbf{x}\right)
\right)
\end{array}
\label{e74}%
\end{equation}
This system is also strongly elliptic.

\begin{theorem}
\label{theorem3}Assume $A(\mathbf{s})$, $\mathbf{s}\in\overline{\Omega}$,
satisfies (\ref{e3}); and without loss of generality, assume
\[
\det J(\mathbf{x})>0,\quad\quad\mathbf{x}\in\overline{B}
\]
Recall $\widetilde{A}\left(  \mathbf{x}\right)  $ as defined by (\ref{e69}).
Then $\widetilde{A}\left(  \mathbf{x}\right)  $ satisfies the strong
ellipticity condition%
\begin{align*}
\xi^{\text{T}}\widetilde{A}(\mathbf{x})\xi &  \geq\widetilde{c}_{0}%
\xi^{\text{T}}\xi,\quad\quad\mathbf{x}\in\overline{B},\quad\xi\in
\mathbb{R}^{d}\\
\widetilde{c}_{0}  &  =c_{0}\lambda_{\ast}\equiv c_{0}\min_{\mathbf{x}%
\in\overline{B}}\lambda_{\min}(\mathbf{x})
\end{align*}
with $\lambda_{\min}(\mathbf{x})$ the smallest eigenvalue of $K(\Phi\left(
\mathbf{x}\right)  )^{\text{T}}K(\Phi\left(  \mathbf{x}\right)  )$ (which
equals the reciprocal of the largest eigenvalue of $J\left(  \mathbf{x}%
\right)  ^{\text{T}}J\left(  \mathbf{x}\right)  $).

\begin{proof}%
\begin{align*}
\xi^{\text{T}}\widetilde{A}(\mathbf{x})\xi &  =\xi^{\text{T}}KAK^{\text{T}}%
\xi=\left(  K^{\text{T}}\xi\right)  ^{\text{T}}A\left(  K^{\text{T}}\xi\right)
\\
&  \geq c_{0}\left(  K^{\text{T}}\xi\right)  ^{\text{T}}\left(  K^{\text{T}%
}\xi\right)  =c_{0}\left\Vert K^{\text{T}}\xi\right\Vert _{2}^{2}%
\end{align*}
In addition,
\[
\left\Vert K(\Phi\left(  \mathbf{x}\right)  )^{\text{T}}\xi\right\Vert
_{2}^{2}\geq\lambda_{\min}(\mathbf{x})\left\Vert \xi\right\Vert _{2}^{2}%
\geq\lambda_{\ast}\left\Vert \xi\right\Vert _{2}^{2}
\]%
\[
\lambda_{\ast}=\min_{\mathbf{x}\in\overline{B}}\lambda_{\min}(\mathbf{x})
\]
with $\lambda_{\min}\left(  \mathbf{x}\right)  $ the smallest eigenvalue of
$K(\Phi\left(  \mathbf{x}\right)  )^{\text{T}}K(\Phi\left(  \mathbf{x}\right)
)$; cf. \cite[p. 488]{atkinson89}. $\left.  {}\right.  $\hfill
\end{proof}
\end{theorem}

\section{Implementation \label{implementation}}

Consider the implementation of the Galerkin method of \S \ref{spec_method} for
the elliptic problem (\ref{e22}) over the unit ball $B$. We are to find the
function $u_{n}\in\mathcal{X}_{n}$ satisfying (\ref{e30}). To do so, we begin
by selecting an orthonormal basis for $\Pi_{n}$, denoting it by $\left\{
\varphi_{1},\dots,\varphi_{N}\right\}  $, with $N\equiv N_{n}=\dim\Pi_{n}$.
Choosing an orthonormal basis is an attempt to have the linear system in
(\ref{e30}) be better conditioned. Next, let
\begin{equation}
\psi_{i}(\mathbf{x})=\left(  1-\left\Vert \mathbf{x}\right\Vert _{2}%
^{2}\right)  \varphi_{i}(\mathbf{x}),\quad\quad i=1,\dots,N_{n} \label{e75}%
\end{equation}
to form a basis for $\mathcal{X}_{n}$.

We seek%
\begin{equation}
u_{n}(\mathbf{x})=\sum_{j=1}^{N}\alpha_{j}\psi_{j}(\mathbf{x}) \label{e76}%
\end{equation}
Then (\ref{e30}) becomes%
\begin{equation}%
\begin{array}
[c]{r}%
{\displaystyle\sum\limits_{k=1}^{N_{n}}}
\alpha_{k}%
{\displaystyle\int_{B}}
\left[
{\displaystyle\sum\limits_{i,j=1}^{d}}
a_{i,j}(\mathbf{x})\dfrac{\partial\psi_{k}(\mathbf{x})}{\partial x_{j}}%
\dfrac{\partial\psi_{\ell}(\mathbf{x})}{\partial x_{i}}+\gamma(\mathbf{x}%
)\psi_{k}(\mathbf{x})\psi_{\ell}(\mathbf{x})\right]  d\mathbf{x}\quad
\quad\quad\medskip\\
=%
{\displaystyle\int_{B}}
f\left(  \mathbf{x}\right)  \psi_{\ell}\left(  \mathbf{x}\right)
\,d\mathbf{x},\quad\quad\ell=1,\dots,N
\end{array}
\label{e78}%
\end{equation}
We need to calculate the orthonormal polynomials and their first partial
derivatives; and we also need to approximate the integrals in the linear
system. For an introduction to the topic of multivariate orthogonal
polynomials, see Dunkl and Xu \cite{DX} and Xu \cite{xu2004}. For multivariate
quadrature over the unit ball in $\mathbb{R}^{d}$, see Stroud \cite{stroud}.

\subsection{The planar case}

The dimension of $\Pi_{n}$ is
\begin{equation}
N_{n}=\frac{1}{2}\left(  n+1\right)  \left(  n+2\right)  \label{e79}%
\end{equation}
For notation, we replace $\mathbf{x}$ with $\left(  x,y\right)  $. How do we
choose the orthonormal basis $\left\{  \varphi_{\ell}(x,y)\right\}  _{\ell
=1}^{N}$ for $\Pi_{n}$? Unlike the situation for the single variable case,
there are many possible orthonormal bases over $B=D$, the unit disk in
$\mathbb{R}^{2}$. We have chosen one that is particularly convenient for our
computations. These are the "ridge polynomials" introduced by Logan and\ Shepp
\cite{Loga} for solving an image reconstruction problem. We summarize here the
results needed for our work.

Let
\[
\mathcal{V}_{n}=\left\{  P\in\Pi_{n}:\left(  P,Q\right)  =0\quad\forall
Q\in\Pi_{n-1}\right\}
\]
the polynomials of degree $n$ that are orthogonal to all elements of
$\Pi_{n-1}$. Then the dimension of $\mathcal{V}_{n}$ is $n+1$; moreover,%
\begin{equation}
\Pi_{n}=\mathcal{V}_{0}\oplus\mathcal{V}_{1}\oplus\cdots\oplus\mathcal{V}_{n}
\label{e100}%
\end{equation}
It is standard to construct orthonormal bases of each $\mathcal{V}_{n}$ and to
then combine them to form an orthonormal basis of $\Pi_{n}$ using the latter
decomposition. \ As an orthonormal basis of $\mathcal{V}_{n}$ we use%
\begin{equation}
\varphi_{n,k}(x,y)=\frac{1}{\sqrt{\pi}}U_{n}\left(  x\cos\left(  kh\right)
+y\sin\left(  kh\right)  \right)  ,\quad\left(  x,y\right)  \in D,\quad
h=\frac{\pi}{n+1} \label{e101}%
\end{equation}
for $k=0,1,\dots,n$. The function $U_{n}$ is the Chebyshev polynomial of the
second kind of degree $n$:%
\begin{equation}
U_{n}(t)=\frac{\sin\left(  n+1\right)  \theta}{\sin\theta},\quad\quad
t=\cos\theta,\quad-1\leq t\leq1,\quad n=0,1,\dots\label{e102}%
\end{equation}
The family $\left\{  \varphi_{n,k}\right\}  _{k=0}^{n}$ is an orthonormal
basis of $\mathcal{V}_{n}$. As a basis of $\Pi_{n}$, we order $\left\{
\varphi_{n,k}\right\}  $ lexicographically based on the ordering in
(\ref{e101}) and (\ref{e100}):%
\[
\left\{  \varphi_{\ell}\right\}  _{\ell=1}^{N}=\left\{  \varphi_{0,0}%
,\,\varphi_{1,0},\,\varphi_{1,1},\,\varphi_{2,0},\,\dots,\,\varphi
_{n,0},\,\dots,\varphi_{n,n}\right\}
\]

Returning to (\ref{e75}), we define%
\begin{equation}
\psi_{n,k}(x,y)=\left(  1-x^{2}-y^{2}\right)  \varphi_{n,k}(x,y) \label{e103}%
\end{equation}
To calculate the first order partial derivatives of $\psi_{n,k}(x,y)$, we need
$U_{n}^{^{\prime}}(t)$. The values of $U_{n}(t)$ and $U_{n}^{^{\prime}}(t)$
are evaluated using the standard triple recursion relations%
\begin{align*}
U_{n+1}(t)  &  =2tU_{n}(t)-U_{n-1}(t)\\
U_{n+1}^{^{\prime}}(t)  &  =2U_{n}(t)+2tU_{n}^{^{\prime}}(t)-U_{n-1}%
^{^{\prime}}(t)
\end{align*}

For the numerical approximation of the integrals in (\ref{e78}), which are
over $B$ being the unit disk, we use the formula%
\begin{equation}
\int_{B}g(x,y)\,dx\,dy\approx\sum_{l=0}^{q}\sum_{m=0}^{2q}g\left(  r_{l}%
,\frac{2\pi\,m}{2q+1}\right)  \omega_{l}\frac{2\pi}{2q+1}r_{l} \label{e106}%
\end{equation}
Here the numbers $\omega_{l}$ are the weights of the $\left(  q+1\right)
$-point Gauss-Legendre quadrature formula on $[0,1]$. Note that
\[
\int_{0}^{1}p(x)dx=\sum_{l=0}^{q}p(r_{l})\omega_{l},
\]
for all single-variable polynomials $p(x)$ with $\deg\left(  p\right)
\leq2q+1 $. The formula (\ref{e106}) uses the trapezoidal rule with $2q+1$
subdivisions for the integration over $\overline{B}$ in the azimuthal
variable. This quadrature is exact for all polynomials $g\in\Pi_{2q}$. This
formula is also the basis of the hyperinterpolation formula discussed in
\cite{hac}.

\subsection{The three-dimensional case}

In $\mathbb{R}^{3}$, the dimension of $\Pi_{n}$ is
\[
N_{n}=\binom{n+3}{3}=\frac{1}{6}\left(  n+1\right)  \left(  n+2\right)
\left(  n+3\right)
\]
Here we choose orthonormal polynomials on the unit ball as described in
\cite{DX},
\begin{align}
\varphi_{m,j,\beta}(\mathbf{x})  &  =c_{m,j}p_{j}^{(0,m-2j+\frac{1}{2}%
)}(2\Vert\mathbf{x}\Vert^{2}-1)S_{\beta,m-2j}\left(  \frac{\mathbf{x}}%
{\Vert\mathbf{x}\Vert}\right)  \medskip\nonumber\\
&  =c_{m,j}\Vert\mathbf{x}\Vert^{m-2j}p_{j}^{(0,m-2j+\frac{1}{2})}%
(2\Vert\mathbf{x}\Vert^{2}-1)S_{\beta,m-2j}\left(  \frac{\mathbf{x}}%
{\Vert\mathbf{x}\Vert}\right)  ,\medskip\label{e1001}\\
j  &  =0,\ldots,\lfloor m/2\rfloor,\quad\beta=0,1,\ldots,2(m-2j),\quad
m=0,1,\ldots,n\nonumber
\end{align}
Here $c_{m,j}=2^{\frac{5}{4}+\frac{m}{2}-j}$ is a constant, and $p_{j}%
^{(0,m-2j+\frac{1}{2})}$, $j\in\mathbb{N}_{0}$, are the normalized Jabobi
polynomials which are orthonormal on $[-1,1]$ with respect to the inner
product
\[
(v,w)=\int_{-1}^{1}(1+t)^{m-2j+\frac{1}{2}}v(t)w(t)\;dt,
\]
see for example \cite{abramowitz}, \cite{gautschi}. The functions
$S_{\beta,m-2j}$ are spherical harmonic functions, and they are given in
spherical coordinates by
\[
S_{\beta,k}(\phi,\theta)=\widetilde{c}_{\beta,k}\left\{
\begin{array}
[c]{ll}%
\cos(\frac{\beta}{2}\phi)T_{k}^{\frac{\beta}{2}}(\cos\theta),\medskip &
\beta\text{ even}\\
\sin(\frac{\beta+1}{2}\phi)T_{k}^{\frac{\beta+1}{2}}(\cos\theta), &
\beta\text{ odd}%
\end{array}
\right.
\]
The constant $\widetilde{c}_{\beta,k}$ is chosen in such a way that the
functions are orthonormal on the unit sphere $S^{2}$ in $\mathbb{R}^{3}$:
\[
\int_{S^{2}}S_{\beta,k}(\mathbf{x})\,S_{\widetilde{\beta},\widetilde{k}%
}(\mathbf{x})\,dS=\delta_{\beta,\widetilde{\beta}}\,\delta_{k,\widetilde{k}}%
\]
The functions $T_{k}^{l}$ are the associated Legendre polynomials, see
\cite{hobson}, \cite{macroberts}. According to (\ref{e75}) we define the basis
for our space of trial functions by
\[
\psi_{m,j,\beta}(\mathbf{x})=(1-\Vert\mathbf{x}\Vert^{2})\varphi_{m,j,\beta
}(\mathbf{x})
\]
and we can order the basis lexicographically. To calculate all of the above
functions we can use recursive algorithms similar to the one used for the
Chebyshev polynomials. These algorithms also allow the calculation of the
derivatives of each of these functions, see \cite{gautschi}, \cite{zhang-jin}

For the numerical approximation of the integrals in (\ref{e78}) we use a
quadrature formula for the unit ball $B$
\begin{align}
\int_{B}g(x)\,dx  &  =\int_{0}^{1}\int_{0}^{2\pi}\int_{0}^{\pi}\widetilde
{g}(r,\theta,\phi)\,r^{2}\sin(\phi)\,d\phi\,d\theta\,dr\approx Q_{q}%
[g]\medskip\nonumber\\
Q_{q}[g]  &  :=\sum_{i=1}^{2q}\sum_{j=1}^{q}\sum_{k=1}^{q}\frac{\pi}%
{q}\,\omega_{j}\,\nu_{k}\widetilde{g}\left(  \frac{\zeta_{k}+1}{2},\frac
{\pi\;i}{2q},\arccos(\xi_{j})\right)  \label{e1004}%
\end{align}
Here $\widetilde{g}(r,\theta,\phi)=g(\mathbf{x})$ is the representation of $g
$ in spherical coordinates. For the $\theta$ integration we use the
trapezoidal rule, because the function is $2\pi-$periodic in $\theta$. For the
$r$ direction we use the transformation
\begin{align*}
\int_{0}^{1}r^{2}v(r)\;dr  &  =\int_{-1}^{1}\left(  \frac{t+1}{2}\right)
^{2}v\left(  \frac{t+1}{2}\right)  \frac{dt}{2}\medskip\\
&  =\frac{1}{8}\int_{-1}^{1}(t+1)^{2}v\left(  \frac{t+1}{2}\right)
\;dt\medskip\\
&  \approx\sum_{k=1}^{q}\underset{_{=:\nu_{k}}}{\underbrace{\frac{1}{8}\nu
_{k}^{\prime}}}v\left(  \frac{\zeta_{k}+1}{2}\right)
\end{align*}
where the $\nu_{k}^{\prime}$ and $\zeta_{k}$ are the weights and the nodes of
the Gauss quadrature with $q$ nodes on $[-1,1]$ with respect to the inner
product
\[
(v,w)=\int_{-1}^{1}(1+t)^{2}v(t)w(t)\,dt
\]
The weights and nodes also depend on $q$ but we omit this index. For the
$\phi$ direction we use the transformation
\begin{align*}
\int_{0}^{\pi}\sin(\phi)v(\phi)\,d\phi &  =\int_{-1}^{1}v(\arccos
(\phi))\,d\phi\medskip\\
&  \approx\sum_{j=1}^{q}\omega_{j}v(\arccos(\xi_{j}))
\end{align*}
where the $\omega_{j}$ and $\xi_{j}$ are the nodes and weights for the
Gauss--Legendre quadrature on $[-1,1]$. For more information on this
quadrature rule on the unit ball in $\mathbb{R}^{3}$, see \cite{stroud}.

Finally we need the gradient in Cartesian coordinates to approximate the
integral in (\ref{e78}), but the function $\varphi_{m,j,\beta}(x)$ in
(\ref{e1001}) is given in spherical coordinates. Here we simply use the chain
rule, with $\mathbf{x}=\left(  x,y,z\right)  $,
\begin{align*}
\frac{\partial}{\partial x}v(r,\theta,\phi)  &  =\frac{\partial}{\partial
r}v(r,\theta,\phi)\cos(\theta)\sin(\phi)-\frac{\partial}{\partial\theta
}v(r,\theta,\phi)\frac{\sin(\theta)}{r\sin(\phi)}\medskip\\
&  +\frac{\partial}{\partial\phi}v(r,\theta,\phi)\frac{\cos(\theta)\cos(\phi
)}{r}%
\end{align*}
and similarly for $\frac{\partial}{\partial y}$ and $\frac{\partial}{\partial
z}$.%

\begin{figure}
[tb]
\begin{center}
\includegraphics[
height=3in,
width=3.9998in
]%
{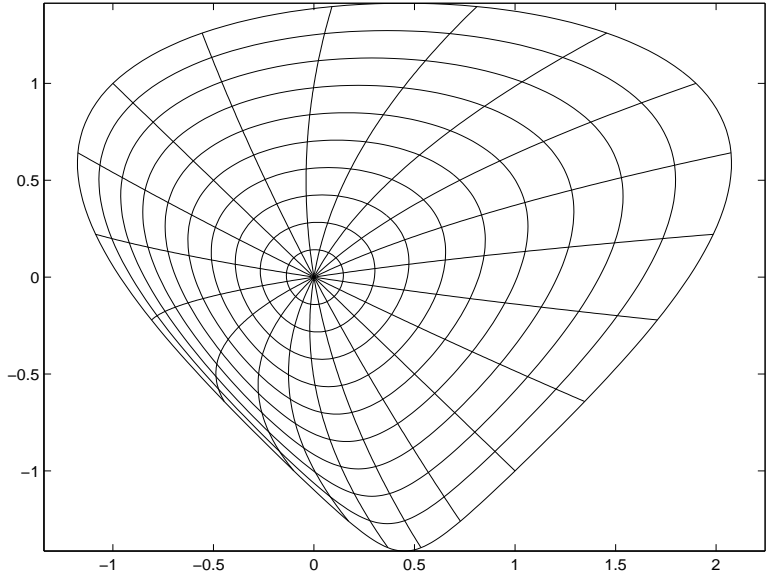}%
\caption{Images of (\ref{e132}), with $a=0.5$, for lines of constant radius
and constant azimuth on the unit disk.}%
\label{fig1}%
\end{center}
\end{figure}
%

\begin{table}[tb] \centering
\caption{Maximum errors in Galerkin solution $u_n$}\label{table1}%
\begin{tabular}
[c]{|c|c|c|c||c|c|c|c|}\hline
$n$ & $N_{n}$ & $\left\Vert u-u_{n}\right\Vert _{\infty}$ & \textit{cond} &
$n$ & $N_{n}$ & $\left\Vert u-u_{n}\right\Vert _{\infty}$ & \textit{cond}%
\\\hline
$2$ & $6$ & $4.41E-1$ & $3.42$ & $14$ & $120$ & $9.95E-6$ & $141.2$\\\hline
$3$ & $10$ & $4.21E-1$ & $4.99$ & $15$ & $136$ & $3.03E-6$ & $165.8$\\\hline
$4$ & $15$ & $1.70E-1$ & $9.27$ & $16$ & $153$ & $8.31E-7$ & $192.8$\\\hline
$5$ & $21$ & $9.63E-2$ & $13.6$ & $17$ & $171$ & $2.09E-7$ & $222.1$\\\hline
$6$ & $28$ & $4.73E-2$ & $20.7$ & $18$ & $190$ & $5.21E-8$ & $253.8$\\\hline
$7$ & $36$ & $1.88E-2$ & $28.5$ & $19$ & $210$ & $1.42E-8$ & $287.9$\\\hline
$8$ & $45$ & $7.24E-3$ & $39.0$ & $20$ & $231$ & $3.53E-9$ & $324.4$\\\hline
$9$ & $55$ & $2.79E-3$ & $50.5$ & $21$ & $253$ & $7.58E-10$ & $363.4$\\\hline
$10$ & $66$ & $9.58E-4$ & $64.7$ & $22$ & $276$ & $1.46E-10$ & $404.9$\\\hline
$11$ & $78$ & $3.20E-4$ & $80.4$ & $23$ & $300$ & $3.36E-11$ & $448.9$\\\hline
$12$ & $91$ & $9.67E-5$ & $98.6$ & $24$ & $325$ & $7.16E-12$ & $495.4$\\\hline
$13$ & $105$ & $3.01E-5$ & $118.7$ & $25$ & $351$ & $1.44E-12$ &
$544.4$\\\hline
\end{tabular}%
\end{table}%

\section{Numerical example\label{num_example}}

Our programs are written in \textsc{Matlab} and can be obtained from the
authors. Our transformations have been so chosen that we can invert explicitly
the mapping $\Phi$, to be able to better construct our test examples. This is
not needed when applying the method; but it simplified the construction of our
test cases. The elliptic equation being solved is
\begin{equation}
Lu(\mathbf{s})\equiv-\Delta u+\gamma(\mathbf{s})u(\mathbf{s})=f(\mathbf{s}%
),\quad\quad\mathbf{s}\in\text{$\Omega$}\subseteq\mathbb{R}^{d} \label{e130}%
\end{equation}
which corresponds to choosing $A=I$. Then we need to calculate%
\begin{equation}%
\begin{array}
[c]{c}%
\widetilde{A}\left(  \mathbf{x}\right)  =K\left(  \Phi\left(  \mathbf{x}%
\right)  \right)  K\left(  \Phi\left(  \mathbf{x}\right)  \right)  ^{\text{T}%
}\smallskip\\
K\left(  \Phi\left(  \mathbf{x}\right)  \right)  =J\left(  \mathbf{x}\right)
^{-1}%
\end{array}
\label{e131}%
\end{equation}

\subsection{The planar case}

For our variables, we replace $\mathbf{x}\in B$ with $\left(  x,y\right)  $,
and we replace $\mathbf{s}\in\Omega$ with $\left(  s,t\right)  $. \ Define the
mapping $\Phi:\overline{B}\rightarrow\overline{\Omega}$ by $\left(
s,t\right)  =\Phi\left(  x,y\right)  $,%
\begin{equation}%
\begin{array}
[c]{l}%
s=x-y+ax^{2}\\
t=x+y
\end{array}
\label{e132}%
\end{equation}
with $0<a<1$. It can be shown that $\Phi$ is a 1-1 mapping from the unit disk
$\overline{B}$. In particular, the inverse mapping $\Psi:\overline{\Omega
}\rightarrow\overline{B}$ is given by%
\begin{equation}%
\begin{array}
[c]{l}%
x=\dfrac{1}{a}\left[  -1+\sqrt{1+a\left(  s+t\right)  }\right]  \medskip\\
y=\dfrac{1}{a}\left[  at-\left(  -1+\sqrt{1+a\left(  s+t\right)  }\right)
\right]
\end{array}
\label{e133}%
\end{equation}
In Figure \ref{fig1}, we give the images in $\overline{\Omega}$ of the circles
$r=j/10$, $j=1,\dots,10$ and the azimuthal lines $\theta=j\pi/10$,
$j=1,\dots,20$.

The following information is needed when implementing the transformation from
$-\Delta u+\gamma u=f$ on $\Omega$ to a new equation on $B$:%
\[
D\Phi=J\left(  x,y\right)  =\left(
\begin{array}
[c]{cc}%
1+2ax & -1\\
1 & 1
\end{array}
\right)
\]%
\[
\det\left(  J\right)  =2\left(  1+ax\right)
\]%
\[
K=\frac{1}{2\left(  1+ax\right)  }\left(
\begin{array}
[c]{cc}%
1 & 1\\
-1 & 1+2ax
\end{array}
\right)
\]%
\[
\widetilde{A}=KK^{T}=\frac{1}{2\left(  1+ax\right)  ^{2}}\left(
\begin{array}
[c]{cc}%
1 & ax\\
ax & 2a^{2}x^{2}+2ax+1
\end{array}
\right)
\]%
\[
\det\left(  J\right)  \widetilde{A}=\frac{1}{1+ax}\left(
\begin{array}
[c]{cc}%
1 & ax\\
ax & 2a^{2}x^{2}+2ax+1
\end{array}
\right)
\]
The latter are the coefficients for the transformed elliptic operator over $B
$, given in (\ref{e68}).

\begin{figure}
[tb]
\begin{center}
\includegraphics[
height=3in,
width=3.9998in
]%
{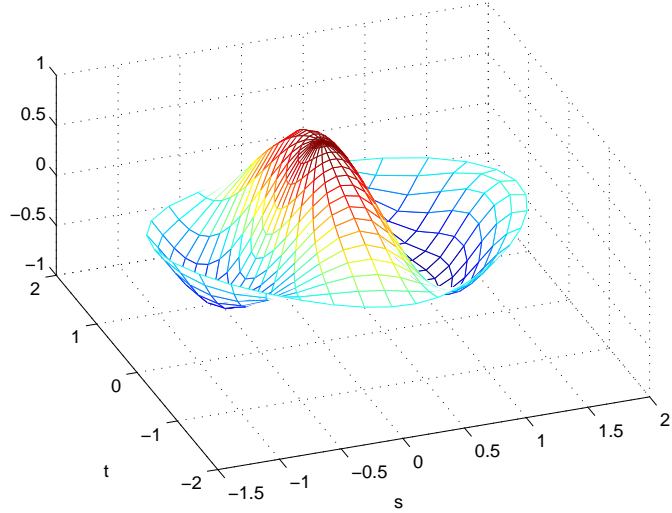}%
\caption{The true solution (\ref{e138})}%
\label{true_soln}%
\end{center}
\end{figure}

We give numerical results for solving the equation%
\begin{equation}
-\Delta u\left(  s,t\right)  +e^{s-t}u\left(  s,t\right)  =f\left(
s,t\right)  ,\quad\quad\left(  s,t\right)  \in\Omega\label{e136}%
\end{equation}
As a test case, we choose%
\begin{equation}
u\left(  s,t\right)  =\left(  1-x^{2}-y^{2}\right)  \cos\left(  \pi s\right)
\label{e138}%
\end{equation}
with $\left(  x,y\right)  $ replaced using (\ref{e133}). The solution is
pictured in Figure \ref{true_soln}. To find $f(s,t)$, we use (\ref{e136}) and
(\ref{e138}). We use the domain parameter $a=0.5$, with $\Omega$ pictured in
Figure \ref{fig1}.%

%

\begin{figure}
[tb]
\begin{center}
\includegraphics[
height=3in,
width=3.9998in
]%
{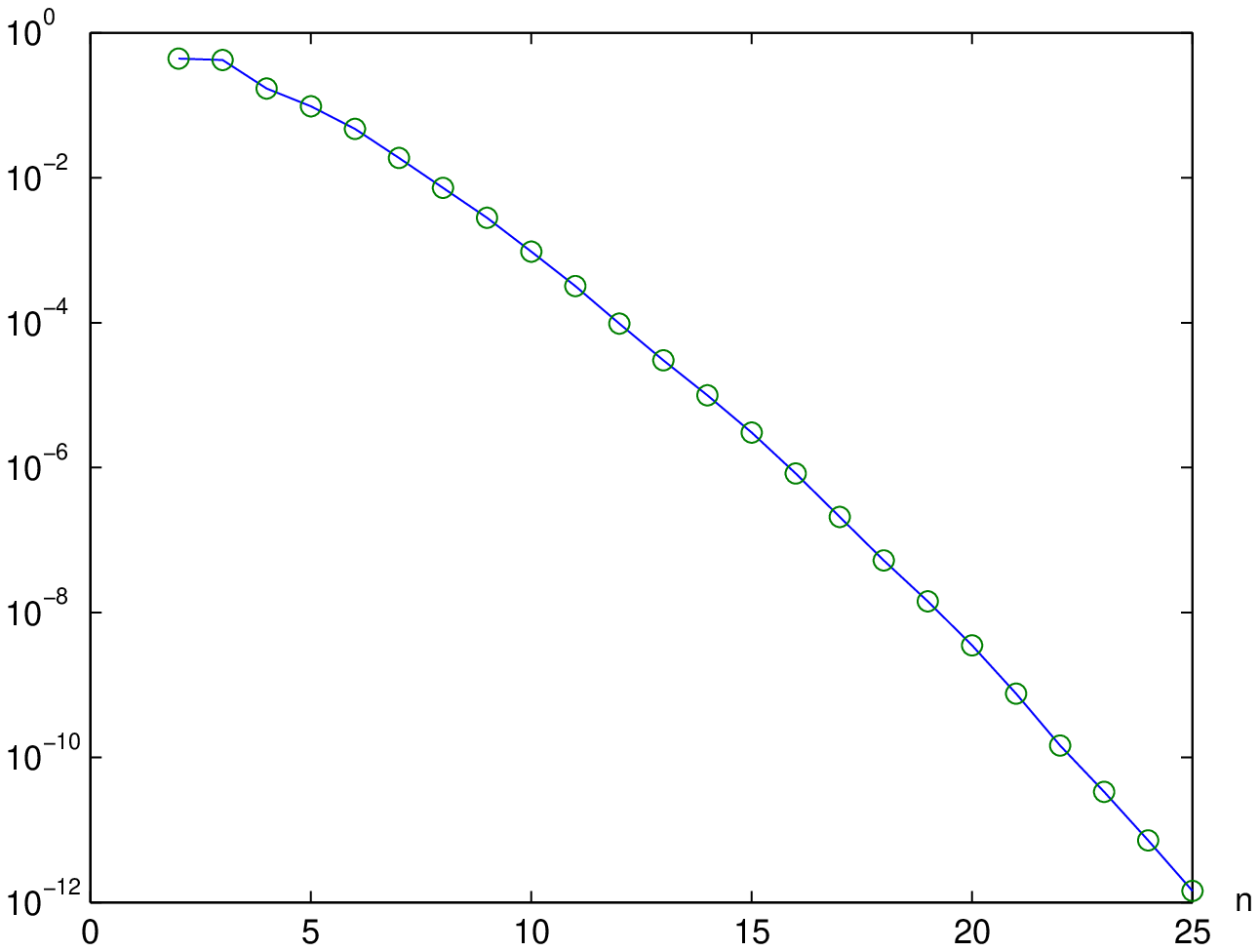}%
\caption{Errors from Table \ref{table1}}%
\label{error_planar}%
\end{center}
\end{figure}
%

\begin{figure}
[tb]
\begin{center}
\includegraphics[
height=3in,
width=3.9998in
]%
{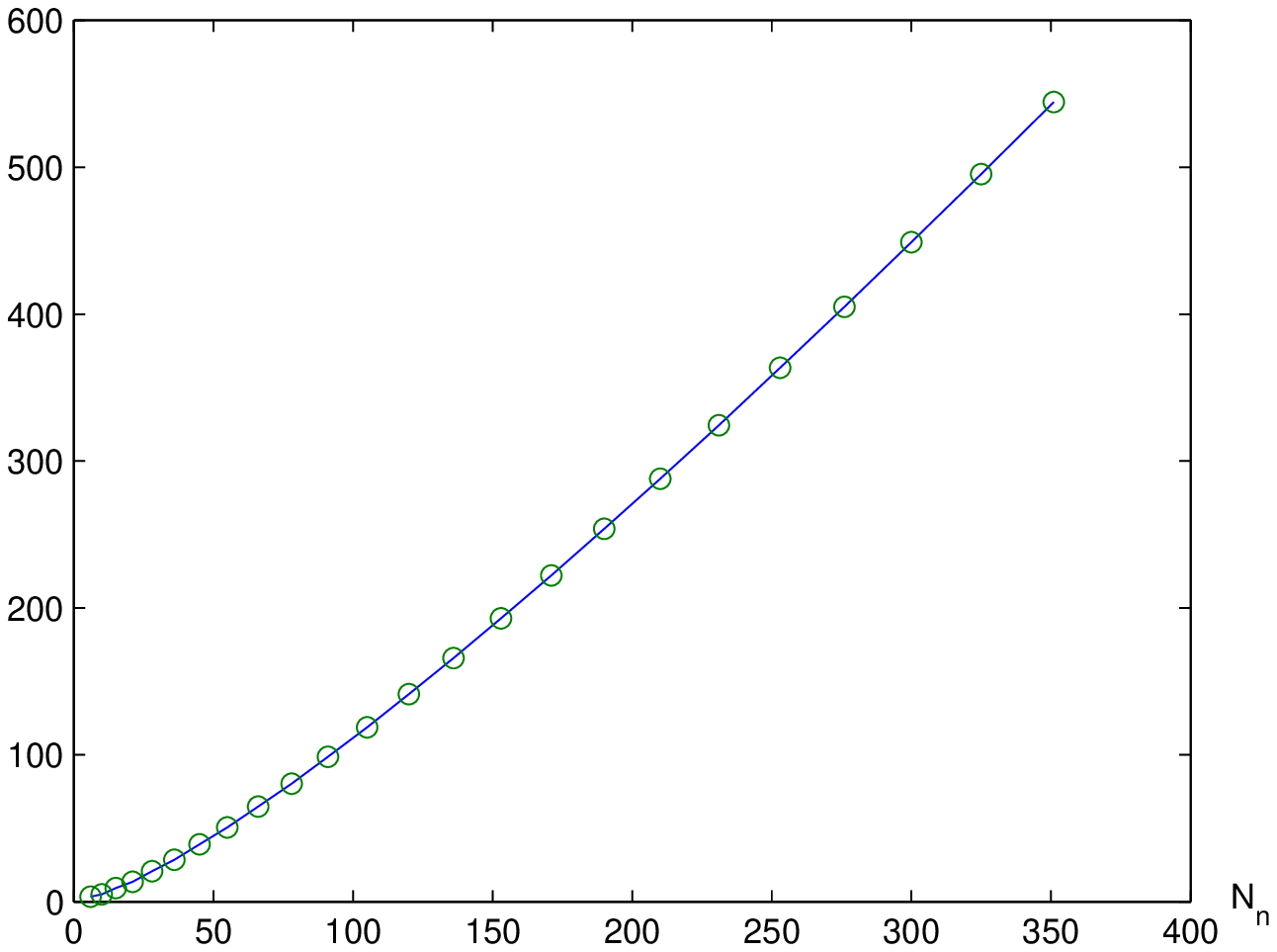}%
\caption{Condition numbers from Table \ref{table1}}%
\label{cond_planar}%
\end{center}
\end{figure}

Numerical results are given in Table \ref{table1}. The integrations in
(\ref{e78}) were performed with (\ref{e106}); and the integration parameter
$q$ ranged from $10$ to $30$. We give the condition numbers of the linear
system (\ref{e78}) as produced in \textsc{Matlab}. To calculate the error, we
evaluate the the numerical solution and the error on the grid
\begin{align*}
\Phi\left(  x_{i,j},y_{i,j}\right)   &  =\Phi\left(  r_{i}\cos\theta_{j}%
,r_{i}\sin\theta_{j}\right)  \smallskip\\
\left(  r_{i},\theta_{j}\right)   &  =\left(  \frac{i}{10},\frac{j\pi}%
{10}\right)  ,\quad\quad i=0,1,\dots10;\quad j=1,\dots20
\end{align*}
The results are shown graphically in Figure \ref{error_planar}. The use of a
semi-log scale demonstrates the exponential convergence of the method as the
degree increases.

To examine experimentally the behaviour of the condition numbers for the
linear system (\ref{e78}), we have graphed the condition numbers from Table
\ref{table1} in Figure \ref{cond_planar}. Note that we are graphing $N_{n}$
vs. the condition number of the associated linear system. \ The graph seems to
indicate that the condition number of the system (\ref{e78}) is directly
proportional to the order of the system, with the order given in (\ref{e79}).%

\begin{figure}[tbp] \centering
\begin{tabular}
[c]{ll}%
{\parbox[b]{2.3333in}{\begin{center}
\includegraphics[
height=1.7504in,
width=2.3333in
]%
{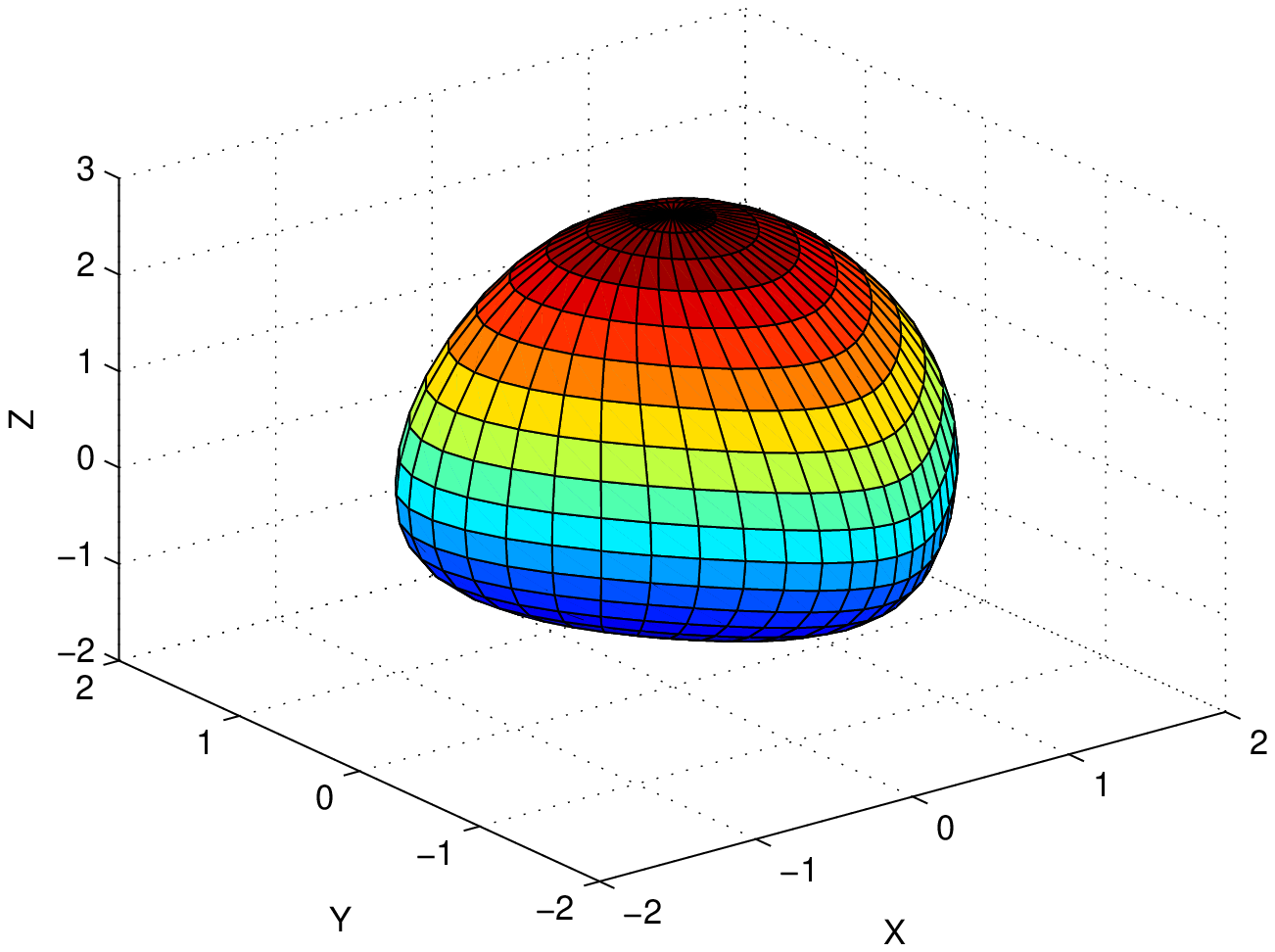}%
\\
{}%
\end{center}}}%
&
{\includegraphics[
height=1.7659in,
width=2.354in
]%
{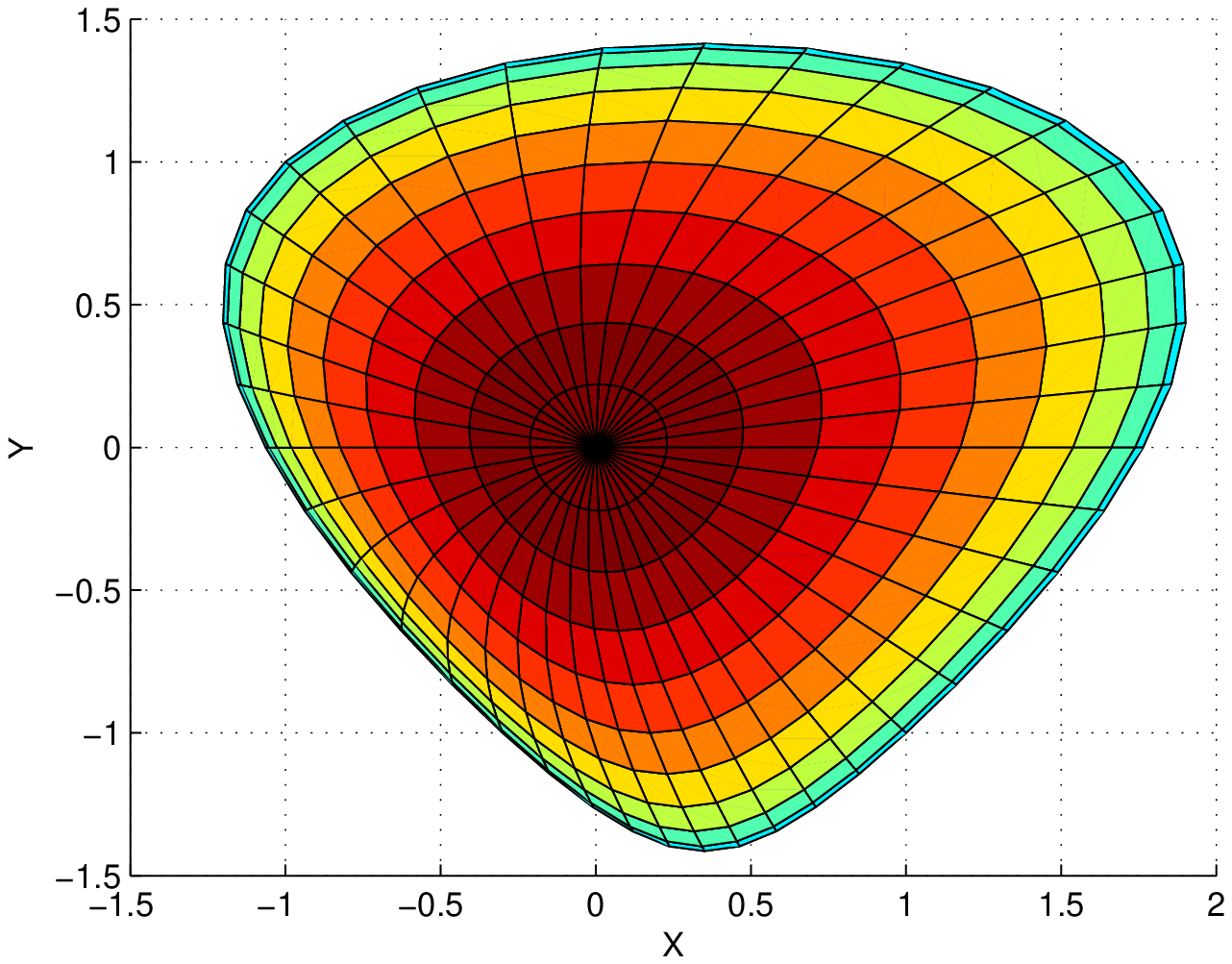}%
}%
\end{tabular}
\caption{Image of (\ref{e1003}) from two different angles, with a=0.7,
b=0.9, for lines of constant $\varphi$ and $\theta$  on the sphere.}\label{fig5}%
\end{figure}%

\subsection{The three-dimensional case}

Here we define the mapping $\Phi:\overline{B}\rightarrow\overline{\Omega}$ by
$(s,t,u)=\Phi(x,y,z)$,
\begin{align}
s  &  =x-y+ax^{2}\nonumber\\
t  &  =x+y\label{e1003}\\
u  &  =2z+bz^{2}\nonumber
\end{align}
$0<a,b<1$, which is an extension of the mapping defined in (\ref{e132}). The
inverse mapping $\Psi:\overline{\Omega}\rightarrow\overline{B}$ is given by
\begin{align*}
x  &  =\dfrac{1}{a}\left[  -1+\sqrt{1+a\left(  s+t\right)  }\right] \\
y  &  =\dfrac{1}{a}\left[  at-\left(  -1+\sqrt{1+a\left(  s+t\right)
}\right)  \right] \\
z  &  =\dfrac{1}{b}\left[  -1+\sqrt{1+bu}\right]
\end{align*}
In Figure \ref{fig5} we show the image of the surface of $\overline{B}$ under
$\Phi$. As in the planar case, we also need
\[
D\Phi(x,y,z):=:J(x,y,z)=\left(
\begin{array}
[c]{ccc}%
1+2ax & -1 & 0\\
1 & 1 & 0\\
0 & 0 & 2+2bz
\end{array}
\right)
\]%
\[
\det(J(x,y,z))=4(1+ax)(1+bz)
\]
and
\begin{align*}
&  \det(J(x,y,z))\widetilde{A}(x,y,z)\\
&  =\det(J(x,y,z)K(x,y,z)K^{T}(x,y,z)\\
&  =4(1+ax)(1+bz)\left(
\begin{array}
[c]{ccc}%
\dfrac{1}{2(1+ax)^{2}} & \dfrac{ax}{2(1+ax)^{2}} & 0\\
\dfrac{ax}{2(1+ax)^{2}} & \dfrac{1+ax+2a^{2}x^{2}}{2(1+ax)^{2}} & 0\\
0 & 0 & \dfrac{1}{4(1+bz)^{2}}%
\end{array}
\right)
\end{align*}
Again, these are the coefficients for the second order term for the
transformed equation on $\overline{B}$, given in (\ref{e68}). We give
numerical results for solving the equation
\[
-\Delta v(s,t,u)+e^{s-t}v(s,t,u)=f(s,t,u),\quad(s,t,u)\in\Omega
\]
and for our test case we choose
\[
v(s,t,u)=\sin\left(  \frac{1}{2}(s-t)\right)  \cdot(1-\Vert\Psi(s,t,u)\Vert
^{2})
\]
where the second term guarantees the Dirichlet boundary conditions on
$\overline{\Omega}$.

Numerical results are given in Table \ref{table2}.%

\begin{table}[tb] \centering
\caption{Maximum errors in Galerkin solution $u_n$}\label{table2}%
\begin{tabular}
[c]{|r|r|r|r|}\hline
$n$ & $N_{n}$ & $\left\Vert u-u_{n}\right\Vert _{\infty}$ & \textit{cond}%
\\\hline
\multicolumn{1}{|c|}{$1$} & \multicolumn{1}{|c|}{$4$} &
\multicolumn{1}{|c|}{$4.98E-1$} & \multicolumn{1}{|c|}{$1.5$}\\\hline
\multicolumn{1}{|c|}{$2$} & \multicolumn{1}{|c|}{$10$} &
\multicolumn{1}{|c|}{$1.99E-1$} & \multicolumn{1}{|c|}{$3.6$}\\\hline
\multicolumn{1}{|c|}{$3$} & \multicolumn{1}{|c|}{$20$} &
\multicolumn{1}{|c|}{$1.78E-1$} & \multicolumn{1}{|c|}{$5.7$}\\\hline
\multicolumn{1}{|c|}{$4$} & \multicolumn{1}{|c|}{$35$} &
\multicolumn{1}{|c|}{$8.22E-2$} & \multicolumn{1}{|c|}{$11.0$}\\\hline
\multicolumn{1}{|c|}{$5$} & \multicolumn{1}{|c|}{$56$} &
\multicolumn{1}{|c|}{$2.18E-2$} & \multicolumn{1}{|c|}{$17.1$}\\\hline
\multicolumn{1}{|c|}{$6$} & \multicolumn{1}{|c|}{$84$} &
\multicolumn{1}{|c|}{$1.34E-2$} & \multicolumn{1}{|c|}{$27.1$}\\\hline
\multicolumn{1}{|c|}{$7$} & \multicolumn{1}{|c|}{$120$} &
\multicolumn{1}{|c|}{$5.95E-3$} & \multicolumn{1}{|c|}{$39.4$}\\\hline
\multicolumn{1}{|c|}{$8$} & \multicolumn{1}{|c|}{$165$} &
\multicolumn{1}{|c|}{$1.60E-3$} & \multicolumn{1}{|c|}{$55.9$}\\\hline
\multicolumn{1}{|c|}{$9$} & \multicolumn{1}{|c|}{$220$} &
\multicolumn{1}{|c|}{$4.85E-4$} & \multicolumn{1}{|c|}{$75.8$}\\\hline
\multicolumn{1}{|c|}{$10$} & \multicolumn{1}{|c|}{$286$} &
\multicolumn{1}{|c|}{$2.56E-4$} & \multicolumn{1}{|c|}{$100.2$}\\\hline
\multicolumn{1}{|c|}{$11$} & \multicolumn{1}{|c|}{$364$} &
\multicolumn{1}{|c|}{$1.44E-4$} & \multicolumn{1}{|c|}{$128.9$}\\\hline
\multicolumn{1}{|c|}{$12$} & \multicolumn{1}{|c|}{$455$} &
\multicolumn{1}{|c|}{$7.85E-5$} & \multicolumn{1}{|c|}{$162.4$}\\\hline
\multicolumn{1}{|c|}{$13$} & \multicolumn{1}{|c|}{$560$} &
\multicolumn{1}{|c|}{$4.19E-5$} & \multicolumn{1}{|c|}{$200.6$}\\\hline
\multicolumn{1}{|c|}{$14$} & \multicolumn{1}{|c|}{$680$} &
\multicolumn{1}{|c|}{$2.33E-5$} & \multicolumn{1}{|c|}{$244.0$}\\\hline
\end{tabular}%
\end{table}%

The integrations in (\ref{e78}) were performed with (\ref{e1004}); and the
integration parameter $q$ was chosen as $q=n+2$. Numerical experiments
indicate that a larger $q$ does not change the results significantly. The
condition numbers for the system (\ref{e78}) were again calculated with
\textsc{Matlab}. An estimation for the error in the maximum norm was
calculated on the grid given by
\[
\left(
\begin{array}
[c]{r}%
x_{i,j,k}\\
y_{i,j,k}\\
z_{i,j,k}%
\end{array}
\right)  =\left(
\begin{array}
[c]{l}%
\frac{i}{21}\sin\left(  \frac{k}{21}\pi\right)  \cos\left(  \frac{2j}{20}%
\pi\right) \\[2mm]%
\frac{i}{21}\sin\left(  \frac{k}{21}\pi\right)  \sin\left(  \frac{2j}{20}%
\pi\right) \\[2mm]%
\frac{i}{21}\cos\left(  \frac{k}{21}\pi\right)
\end{array}
\right)  ,\quad i,k=1,\ldots,20,\quad j=1,\ldots,40.
\]
The error for the Galerkin method is shown in Figure \ref{error_spatial} and
the development of the condition number is shown in Figure \ref{cond_spatial}.
Again the numerical experiment seems to indicate an exponential convergence of
the method and a linear growth of the condition numbers with respect to the
number of degrees of freedom $N_{n}$ of the linear system (\ref{e78}).

\begin{center}%
\begin{figure}
[tb]
\begin{center}
\includegraphics[
height=3in,
width=3.9998in
]%
{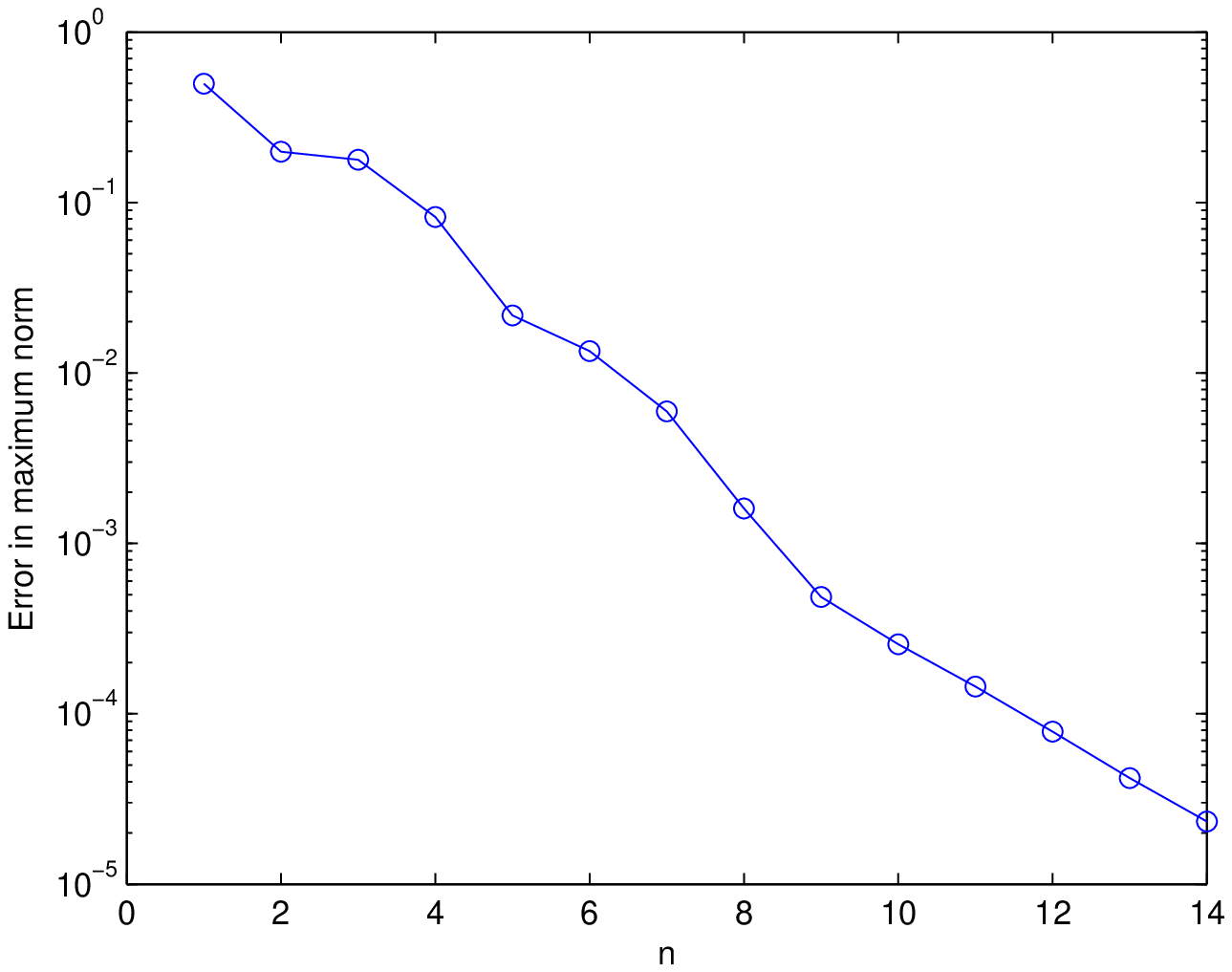}%
\caption{Errors from Table \ref{table2}}%
\label{error_spatial}%
\end{center}
\end{figure}
\begin{figure}
[tb]
\begin{center}
\includegraphics[
height=3in,
width=3.9998in
]%
{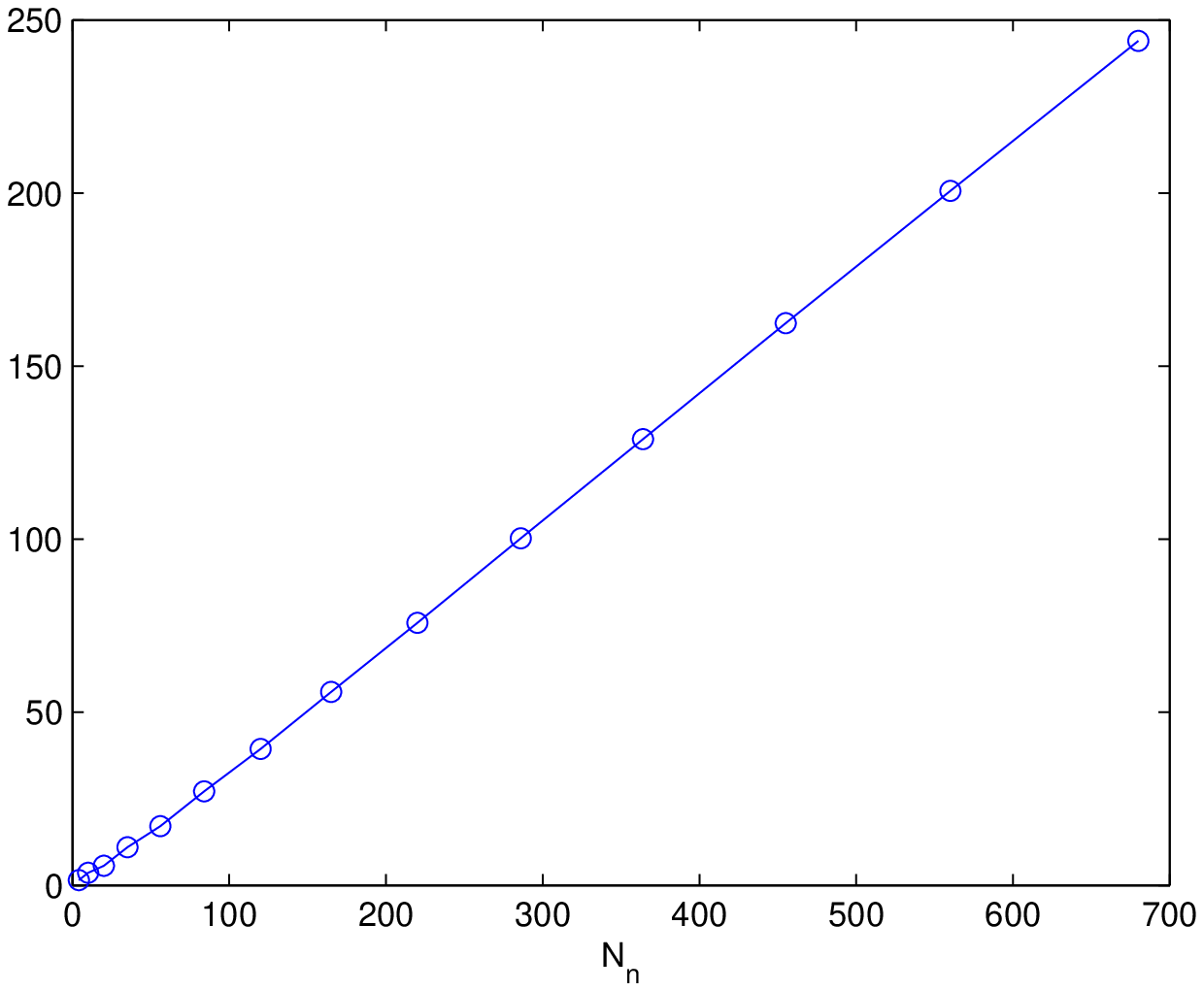}%
\caption{Condition numbers from Table \ref{table2}}%
\label{cond_spatial}%
\end{center}
\end{figure}

\end{center}

\noindent\textsc{Additional Remarks}. We present and study a spectral method
for the Neumann problem%
\begin{align*}
-\Delta u(\mathbf{s})+\gamma(\mathbf{s})u(\mathbf{s})  &  =f(\mathbf{s}%
),\quad\quad\mathbf{s}\in\text{$\Omega$}\subseteq\mathbb{R}^{d}\medskip\\
\frac{\partial u(\mathbf{s})}{\partial\mathbf{n}_{\mathbf{s}}}  &
=g(\mathbf{s}),\quad\quad\mathbf{s}\in\partial\text{$\Omega$}%
\end{align*}
in a forthcoming paper. \ We are also investigating the behaviour of the
condition number for the linear system (\ref{e78}) associated with our
spectral method, attempting to prove that it has size $\mathcal{O}(N_{n})$,
consistent with the numbers shown in Tables \ref{table1} and \ref{table2}.

Our earlier numerical examples use given
\[
\Phi:\overline{B}_{d}\underset{onto}{\overset{1-1}{\longrightarrow}}%
\overline{\Omega}%
\]
chosen to be nontrivial and illustrative. In general, however, when given a
smooth mapping%
\[
\varphi:S\underset{onto}{\overset{1-1}{\longrightarrow}}\partial\Omega,
\]
it may not be clear as to how to extend $\varphi$ to $\Phi$ over $B$. In some
cases there is an obvious choice, as when $\Omega$ is an ellipsoid. We are
investigating schemes to produce continuously differentiable extensions $\Phi$
which satisfy%
\[
\min_{\mathbf{x}\in\overline{B}}\,\left\vert \det J(\mathbf{x})\right\vert >0
\]
and for which $J(\mathbf{x})$ is easily computable.

\end{document}